\documentclass[a4paper, 12pt]{article}
\usepackage{}

\usepackage{mathrsfs}
\usepackage{epsfig}
\usepackage{amsmath}
\usepackage{amssymb}
\usepackage{latexsym}
\usepackage{amsfonts}
\usepackage{amsthm}

\usepackage[numbers,sort&compress]{natbib}
\usepackage{graphicx}
\ifx\pdfoutput\undefined  \else
 \fi

\usepackage[centerlast]{caption2}

\usepackage{color}

\usepackage{txfonts}
\usepackage{amssymb}
\usepackage{mathrsfs}
\usepackage{amsfonts}

\newtheorem{theorem}{Theorem}[section]
\newtheorem{lemma}[theorem]{Lemma}

\newtheorem{conj}[theorem]{Conjecture}

\usepackage{latexsym,amssymb}
\usepackage{amsmath}

\newcommand{\F}{{\mathbb F}}

\begin{document}

\title{Davenport constant of the multiplicative semigroup of the ring $\mathbb{Z}_{n_1}\oplus\cdots\oplus \mathbb{Z}_{n_r}$}

\author{
Guoqing Wang$^a$\thanks{Email: gqwang1979@aliyun.com} \ \ \ \ \ \
Weidong Gao$^b$\thanks{Email :
wdgao1963@aliyun.com}\\
{\small $^a$Department of Mathematics, Tianjin Polytechnic University, Tianjin, 300387, P. R. China}\\
{\small $^b$Center for Combinatorics, LPMC-TJKLC, Nankai University, Tianjin 300071, P. R. China}\\
}

\date{}
\maketitle

\begin{abstract} Given a finite commutative semigroup $\mathcal{S}$ (written additively), denoted by ${\rm D}(\mathcal{S})$ the Davenport constant of $\mathcal{S}$,
namely the least positive integer $\ell$ such that for any $\ell$ elements $s_1,\ldots,s_{\ell}\in \mathcal{S}$ there exists a set $I\subsetneq [1,\ell]$ for which $\sum_{i\in I} s_i=\sum_{i=1}^{\ell} s_i$.

Then, for any integers $r\geq 1, n_1,\ldots,n_r>1$, let
$R=\mathbb{Z}_{n_1}\oplus\cdots\oplus \mathbb{Z}_{n_r}$ be the direct sum of these $r$ residue class rings $\mathbb{Z}_{n_1}, \ldots,\mathbb{Z}_{n_r}$.
Moreover, let $\mathcal{S}_R$ be the multiplicative semigroup of the ring $R$, and ${\rm U}(\mathcal{S}_R)$ the group of units of $\mathcal{S}_R$.
In this paper, we prove that $${\rm D}({\rm U}(\mathcal{S}_R))+P_2\leq {\rm D}(\mathcal{S}_R)\leq {\rm D}({\rm U}(\mathcal{S}_R))+\delta,$$ where $P_2=\sharp\{i\in [1,r]: 2 \parallel n_i\}$ and $\delta=\sharp\{i\in [1,r]: 2\mid n_i\}.$
This corrects our previous published wrong result on this problem.

\end{abstract}

\noindent{\sl Key Words}: Davenport constant;  Multiplicative semigroups; Residue class rings

\section {Introduction}

The Davenport constant of any finite abelian group
$G$, denoted ${\rm D}(G)$, is the smallest positive integer $\ell$  such that,
every sequence $T$ of elements in $G$ of length at least $\ell$ contains a nonempty
subsequence $T'$ with the sum of all terms of $T'$ equaling the identity element of $G$. Though
attributed to H. Davenport who proposed the study
of this constant in 1965, K. Rogers \cite{rog1} in 1963 pioneered the investigation of a combinatorial invariant associated with an arbitrary finite abelian group $G$.
The Davenport constant is a central concept of zero-sum theory and has been investigated by many researchers in the scope of finite abelian groups.

In 2008, the two authors  of this manuscript formulated the definition of the Davenport constant of finite commutative semigroups. Subsequently, some related additive results are obtained in the setting of semigroups (see \cite{AdhikariGaoWang14, wangDavenportII, wangErdosBugess, wangAddtiveirreducible, wang-zhang-wang-qu, wang-zhang-qu}).

\noindent \textbf{Definition A.} (\cite{wanggao}) \ {\sl Let $\mathcal{S}$ be a finite commutative semigroup. Let $T$ be a sequence of terms from the semigroup $\mathcal{S}$. We call $T$  reducible if $T$ contains a proper subsequence $T'$ ($T'\neq T$) such that the sum of all terms of $T'$ equals the sum of all terms of $T$. Define the Davenport constant of the semigroup $\mathcal{S}$, denoted ${\rm D}(\mathcal{S})$, to be the smallest $\ell\in \mathbb{N}$ such that every sequence $T$ of length at least $\ell$ of elements in $\mathcal{S}$ is reducible.}

In 2006, A. Geroldinger and F. Halter-Koch had introduced another combinatorial
invariant, which they denoted by ${\rm d}(\mathcal{S})$, (see Definition 2.8.12 in \cite{GH}), and now is called the small
Davenport constant of $\mathcal{S}$ after \cite{AdhikariGaoWang14}; this is closely related to the Davenport constant of $\mathcal{S}$, as it
is known from Proposition 1.2 in \cite{AdhikariGaoWang14} that ${\rm D}(\mathcal{S})={\rm d}(\mathcal{S})+1$ for any finite commutative semigroup $\mathcal{S}$.

It is embarrass that the following obtained result on the Davenport constant of finite commutative semigroups was observed to be incorrect when the ring $R$ is even in some cases.

\noindent \textbf{Theorem B.} (\cite{wanggao}) \ {\sl For integers $r\geq 1$, $n_1,\ldots,n_r>1$, let $R=\mathbb{Z}_{n_1}\oplus\cdots\oplus \mathbb{Z}_{n_r}$. Let $\mathcal{S}_R$ be the multiplicative semigroup of the ring $R$ and ${\rm U}(\mathcal{S}_R)$ be the group of units of $\mathcal{S}_R$. Then $${\rm D}(\mathcal{S}_R)={\rm D}({\rm U}(\mathcal{S}_R))+P_2,$$ where $P_2=\sharp\{i\in [1,r]: 2 \parallel n_i\}$.}

It is high time to correct this mistake.
In this paper, we shall prove the following result by employing a different method with the previous one used in \cite{wanggao}.

\begin{theorem}\label{Theorem main} \ For integers $r\geq 1$, $n_1,\ldots,n_r>1$, let $R=\mathbb{Z}_{n_1}\oplus\cdots\oplus \mathbb{Z}_{n_r}$. Let $\mathcal{S}_R$ be the multiplicative semigroup of the ring $R$ and ${\rm U}(\mathcal{S}_R)$ be the group of units of $\mathcal{S}_R$. Then $${\rm D}({\rm U}(\mathcal{S}_R))+P_2\leq {\rm D}(\mathcal{S}_R)\leq {\rm D}({\rm U}(\mathcal{S}_R))+\delta,$$ where $P_2=\sharp\{i\in [1,r]: 2 \parallel n_i\}$ and $\delta=\sharp\{i\in [1,r]: 2\mid n_i\}.$
\end{theorem}

\section{The preliminaries}

\noindent $\bullet$ In the rest of this manuscript, we shall always admit that $\mathcal{S}$ is a {\sl unitary} finite commutative semigroup.

The operation on $\mathcal{S}$ is denoted by $+$.
The identity element of $\mathcal{S}$, denoted $0_{\mathcal{S}}$, is the unique element $e$ of
$\mathcal{S}$ such that $e+a=a$ for every $a\in \mathcal{S}$. Let
${\rm U}(\mathcal{S})=\{a\in \mathcal{S}: a+a'=0_{\mathcal{S}} \mbox{ for some }a'\in \mathcal{S}\}$ be the group of units
of $\mathcal{S}$. For any element $c\in\mathcal{S}$, let $${\rm St}(c)=\{a\in {\rm U}(\mathcal{S}): a+c=c\}$$ denote the stabilizer of $c$ in the group ${\rm U}(\mathcal{S})$.
Green's preorder on the semigroup $\mathcal{S}$, denoted $\leqq_{\mathcal{H}}$, is defined by
$$a \leqq_{\mathcal{H}} b\Leftrightarrow a=b \ \ \mbox{or}\ \ a=b+c \mbox{ for some }c\in \mathcal{S}.$$ Green's congruence on $\mathcal{S}$, denoted
$\mathcal{H}$, is defined by:
$$a \ \mathcal{H} \ b \Leftrightarrow a \ \leqq_{\mathcal{H}} \ b \mbox{ and } b \ \leqq_{\mathcal{H}} \ a.$$
We write $a<_{\mathcal{H}} b$ to mean that $a \leqq_{\mathcal{H}} b$ but $a \ \mathcal{H} \ b$ does not hold.

A sequence $T$ of $\mathcal{S}$ is denoted by $T=a_1a_2\cdots a_{\ell}=\coprod\limits_{a\in \mathcal{S}}a^{[{\rm v}_a(T)]},$ where $[{\rm v}_a(T)]$ means that the element $a$ occurs ${\rm v}_a(T)$ times in the sequence $T$.
By $\cdot$ we denote the operation to join sequences.
By $|T|$ we denote the length of the sequence, i.e., $|T|=\sum\limits_{a\in \mathcal{S}}{\rm v}_a(T)=\ell.$
Let $T_1,T_2$ be two sequences of $\mathcal{S}$. We call $T_2$
a subsequence of $T_1$ if ${\rm v}_a(T_2)\leq {\rm v}_a(T_1)$ for every element $a\in \mathcal{S}$, denoted by $T_2\mid T_1.$ In particular, if $T_2\neq T_1$, we call $T_2$ a {\sl proper} subsequence of $T_1$, and write $T_3=T_1  T_2^{[-1]}$ to mean the unique subsequence of $T_1$ with $T_2\cdot T_3=T_1$. Let $\varepsilon$ be the {\sl empty sequence}. In particular, the empty sequence $\varepsilon$ is a proper subsequence of any nonempty sequence.
If $T$ is a nonempty sequence, then we let $\sigma(T)=\sum\limits_{a\in\mathcal{S}} [v_a(T)] a.$
We also define $\sigma(\varepsilon)=0_\mathcal{S}.$
We say that $T$ is {\it
reducible} if $\sigma(T')=\sigma(T)$ for some proper subsequence $T'$ of $T$.

In what follows, we  denote $$R=\mathbb{Z}_{n_1}\oplus\cdots\oplus \mathbb{Z}_{n_r}$$
where $r\geq 1, n_1,\ldots, n_r>1$ are integers, and denote $\mathcal{S}_R$ to be the multiplicative semigroup of the ring $R$. For any element ${\bf a}\in \mathcal{S}_{R}$, we denote $\theta_{\bf a}=(a_1,\ldots,a_r)\in [1,n_1]\times \cdots\times [1,n_r]$ be the unique $r$-tuple of integers such that $(\overline{a}_1,\ldots,\overline{a}_r)$ is the corresponding form of the element ${\bf a}$ in the ring $R$. Let $\kappa_i(\theta_{\bf a})=a_i$ where $i\in [1,r]$.
We remark that, since the operation of the semigroup $\mathcal{S}_R$ is always denoted by $+$,
for any elements ${\bf a}, {\bf b}$ and ${\bf c}$ of $\mathcal{S}_R$, ${\bf a}+{\bf b}={\bf c}$ holds in $\mathcal{S}_R$ if and only if $\kappa_i(\theta_a)*\kappa_i(\theta_b)\equiv \kappa_i(\theta_c)\pmod {n_i} \ \ \mbox{ for all } i=1,2,\ldots,r.$

Here, the following two lemmas are necessary.

\begin{lemma}(\cite{GH}, Lemma 6.1.3) \label{Lemma recusive Davenport constant} \ Let $G$ be a finite abelian group, and let $H$ be a subgroup of $G$. Then, ${\rm D}(G)\geq {\rm D}(G/H)+{\rm D}(H)-1$.
\end{lemma}

\noindent $\bullet$ For any prime $p$ and any integer $n\neq 0$, let ${\rm pot}_p(n)$ be largest integer $k$ such that $p^k$ divides $n$.

\begin{lemma}\label{Lemma three conclusions} \ Let $a$ and $b$ be two elements of $\mathcal{S}_R$. Then the following conclusions hold:

{\rm (i)} \ If $a \leqq_{\mathcal{H}} b$, then $\gcd(\kappa_i(\theta_b),n_i)\mid \gcd(\kappa_i(\theta_a),n_i)$ for each $i\in [1,r]$, and ${\rm St}(b)\subseteq{\rm St}(a)$;

{\rm (ii)} \ $a \ {\mathcal{H}} \ b$ if and only if $\gcd(\kappa_i(\theta_b),n_i)= \gcd(\kappa_i(\theta_a),n_i)$ for each $i\in [1,r]$;

{\rm (iii)} \ Suppose $a <_{\mathcal{H}} b$. If there exists some index $t\in [1,r]$ such that $${\rm pot}_{p}(\gcd(\kappa_t(\theta_b), \ n_t))<{\rm pot}_{p}(\gcd(\kappa_t(\theta_a), \ n_t))$$ for some prime $p>2$, or $${\rm pot}_2(\gcd(\kappa_t(\theta_b), \ n_t))<{\rm pot}_2(\gcd(\kappa_t(\theta_a), \ n_t))< {\rm pot}_2(n_t),$$
then ${\rm St}(b)\subsetneq {\rm St}(a)$.
\end{lemma}

\medskip

\begin{proof}
(i) \   Note that $a \leqq_{\mathcal{H}} b$ implies $a=b+c$ for some $c\in \mathcal{S}_R$ since $\mathcal{S}_R$ is unitary, and equivalently, $\kappa_i(\theta_a)\equiv \kappa_i(\theta_b)*\kappa_i(\theta_c)\pmod {n_i}$  for each $i\in [1,r].$ Then  Conclusion (i) follows from a routine verification.

(ii) \  Assume $\gcd(\kappa_i(\theta_b),n_i)= \gcd(\kappa_i(\theta_a),n_i)$ for each $i\in [1,r]$.
It follows that there exist integers $c_i, c_i'\in [0,n_i-1]$ such that $\kappa_i(\theta_a)* c_i\equiv \kappa_i(\theta_b)\pmod {n_i}$
and
$\kappa_i(\theta_b)* c_i'\equiv \kappa_i(\theta_a)\pmod {n_i},$ where $i\in [1,r]$.
Take elements $c, c'\in \mathcal{S}_R$ such that $\kappa_i(\theta_c)=c_i$ and $\kappa_i(\theta_{c'})=c_i'$ for each $i\in [1,r]$. It follows that $a+c=b$ and $b+c'=a$, and so $a \ \mathcal{H} \ b.$
While, the converse follows from Conclusion {\rm (i)}  immediately.

(iii) \ Let $q$ be the largest prime with
\begin{equation}\label{equation pot p a> pot p b}
{\rm pot}_{q}(\gcd(\kappa_t(\theta_b), \ n_t))<{\rm pot}_{q}(\gcd(\kappa_t(\theta_a), \ n_t)).
\end{equation}
Let \begin{equation}\label{equation alpha=}
\alpha={\rm pot}_{q}(\gcd(\kappa_t(\theta_a), \ n_t)).
\end{equation}
Take an element $d\in \mathcal{S}_R$ with
\begin{equation}\label{equation kappi(d)=1 for almost all}
\kappa_i(\theta_d)=1 \ \mbox{ for each } \ i\in [1,r]\setminus \{t\},
\end{equation} and with
\begin{equation}\label{equation kt(thetad)}
\kappa_t(\theta_d)\equiv\left\{ \begin{array}{ll}
2\frac{n_t}{q^{\alpha}}+1 \pmod {n_t} & \textrm{if $q>2$ and $\gcd(2\frac{n_t}{q^{\alpha}}+1, n_t)=1$;}\\
\frac{n_t}{q^{\alpha}}+1 \pmod {n_t} & \textrm{if otherwise.}\\
\end{array} \right.
\end{equation}

Then we have the following.

\noindent \textbf{Assertion A.} \  $d\in {\rm U}(\mathcal{S}_R)$.

{\sl Proof of Assertion A.} \ By \eqref{equation kappi(d)=1 for almost all}, it suffices to establish that $\gcd(\kappa_t(\theta_d), n_t)=1$. If $q=2$, then $\kappa_t(\theta_d)\equiv \frac{n_t}{q^{\alpha}}+1 \pmod {n_t}$, and so $\gcd(\kappa_t(\theta_d), n_t)=1$ follows from the hypothesis, ${\rm pot}_2(\gcd(\kappa_t(\theta_b), \ n_t))<{\rm pot}_2(\gcd(\kappa_t(\theta_a), \ n_t))< {\rm pot}_2(n_t),$ immediately. Assume $$q>2.$$ We show that  $\gcd(\frac{n_t}{q^{\alpha}}+1,n_t)=1$
or $\gcd(2\frac{n_t}{q^{\alpha}}+1, n_t)=1$. Suppose to the contrary that $\gcd(\frac{n_t}{q^{\alpha}}+1,n_t)>1$ and $\gcd(2\frac{n_t}{q^{\alpha}}+1, n_t)>1$.
Since $w\not\mid \frac{n_t}{q^{\alpha}}+1$ and $w\not\mid 2\frac{n_t}{q^{\alpha}}+1$ for every prime divisor $w$ of $n_t$ with $w\neq q$,
it follows that $q\mid \frac{n_t}{q^{\alpha}}+1$ and $q\mid 2\frac{n_t}{q^{\alpha}}+1$, which implies that
$q\mid 2(\frac{n_t}{q^{\alpha}}+1)-(2\frac{n_t}{q^{\alpha}}+1)=1$, which is absurd. This proves Assertion A.
\qed

By \eqref{equation pot p a> pot p b}, \eqref{equation alpha=},  \eqref{equation kappi(d)=1 for almost all} and \eqref{equation kt(thetad)}, we check that $\kappa_i(\theta_d)*\kappa_i(\theta_a)\equiv \kappa_i(\theta_a)\pmod {n_i}\mbox{ for each } \ i\in [1,r],$
and $\kappa_t(\theta_d)*\kappa_t(\theta_b)\not\equiv \kappa_t(\theta_b)\pmod {n_t}.$ That is, $d+a=a$ and $d+b\neq b.$ Combined with Assertion A, we have that $d\in {\rm St}(a)\setminus {\rm St}(b)$ proved. This completes the proof of Lemma \ref{Lemma three conclusions}.  \end{proof}

\section{The proof of Theorem \ref{Theorem main}}

\noindent {\sl Proof of Theorem  \ref{Theorem main}.} \
We first prove that ${\rm D}(\mathcal{S}_R)\geq {\rm D}({\rm U}(\mathcal{S}_R))+P_2$. Assume without loss of generality that $$\{i\in [1,r]: 2\parallel n_i\}=[1,P_2].$$
Take an irreducible sequence $A$ of terms from ${\rm U}(\mathcal{S}_R)$ of length ${\rm D}({\rm U}(\mathcal{S}_R))-1$. Let
\begin{equation}\label{equation B=Acdot}
B=A\cdot \coprod\limits_{i=1}^{P_2} b_i
\end{equation}
 where $\theta_{b_i}=(1,\ldots,1,2,1\ldots,1)$ with $2$ appears at the $i$-th location for each $i\in [1,P_2]$. Now we show that $B$ is an irreducible sequence. Suppose to the contrary that $B$ contains a {\sl proper} subsequence $B'$ with
\begin{equation}\label{equation sigma(B')=sigma(B)}
\sigma(B')=\sigma(B).
\end{equation}
Since $\sigma(B') \ \mathcal{H}\ \sigma(B)$,  it follows from Lemma \ref{Lemma three conclusions} (ii) that $\coprod\limits_{i=1}^{P_2} b_i \mid B'$, say
\begin{equation}\label{equation A'=}
B'=A'\cdot \coprod\limits_{i=1}^{P_2} b_i
\end{equation}
 where $A'$ is a proper subsequence of $A$. By \eqref{equation B=Acdot}, \eqref{equation sigma(B')=sigma(B)} and \eqref{equation A'=}, we derive that
 \begin{equation}\label{equation 2ki=2ki}
 2\kappa_i(\theta_{\sigma(A')})\equiv 2\kappa_i(\theta_{\sigma(A)})\pmod {n_i}\ \mbox{ for each }\ i\in [1,P_2]
 \end{equation}
and
 \begin{equation}\label{equation kj=kj}\kappa_j(\theta_{\sigma(A')})\equiv \kappa_j(\theta_{\sigma(A)})\pmod {n_j}\ \mbox{ for each }\ j\in [P_2+1,r].
 \end{equation} Since $\kappa_i(\theta_{\sigma(A')})$ and $\kappa_i(\theta_{\sigma(A)})$ are odd, it follows from \eqref{equation 2ki=2ki} that  $\kappa_i(\theta_{\sigma(A')})\equiv \kappa_i(\theta_{\sigma(A)})\pmod {n_i}$, where $i\in [1,P_2].$ Combined with \eqref{equation kj=kj}, we have that $\sigma(A')=\sigma(A)$, which is a contradiction.  This proves that $B$ is irreducible, and so ${\rm D}(\mathcal{S}_R)\geq |B|+1={\rm D}({\rm U}(\mathcal{S}_R))+P_2$.
Now it remains to show that ${\rm D}(\mathcal{S}_R)\leq {\rm D}({\rm U}(\mathcal{S}_R))+\delta$.

Let $T=a_1\cdot a_2\cdots a_{\ell}$ be an arbitrary sequence of term from the semigroup $\mathcal{S}_R$ of length $\ell={\rm D}({\rm U}(\mathcal{S}_R))+\delta.$ It suffices to show that $T$ contains a {\sl proper} subsequence $T'$ with $\sigma(T')=\sigma(T)$.
Take a shortest subsequence $V$ of $T$ such that
\begin{equation}\label{equation sigma(V)Hsigma(T)}
\sigma(V) \ \mathcal{H} \ \sigma(T).
\end{equation}
We may assume without loss of generality that $$V=a_1 \cdot a_2\cdots a_t\ \ \ \mbox{where} \ \ t\in [0,\ell].$$ If $t=0$, i.e., $V=\varepsilon$, then $\sigma(V)=0_{\mathcal{S}_R}$, which implies that $\sigma(T)\in {\rm U}(\mathcal{S}_R)$ by \eqref{equation sigma(V)Hsigma(T)}. It follows that $T$ is a sequence of terms from the group ${\rm U}(\mathcal{S}_R)$ and of length $|T|={\rm D}({\rm U}(\mathcal{S}_R))+\delta\geq {\rm D}({\rm U}(\mathcal{S}_R))$, and thus, $T$ is reducible, we are done. Hence, we assume that $$t>0.$$
By the minimality of $|V|$, we derive that
\begin{equation}\label{equation green's chain}
0_{\mathcal{S}_R}>_{\mathcal{H}}a_1>_{\mathcal{H}}(a_1+a_2)>_{\mathcal{H}}\cdots>_{\mathcal{H}}\sum_{i=1}^t a_i.
\end{equation}
Recall that an empty sum of elements of $\mathcal{S}_R$ is taken equal to $0_{\mathcal{S}_R}$.
Denote
$K_i={\rm St}(\sum\limits_{j=1}^i a_j)$ where $i\in [0,t].$
Note that $K_i$ is a subgroup of ${\rm U}(\mathcal{S}_R)$ for each $i\in [0,t]$. Combined with
\eqref{equation green's chain} and Lemma \ref{Lemma three conclusions} (i), we have that
\begin{equation}\label{equation subgroup chains}
K_0\subseteq K_1\subseteq K_2\subseteq \cdots\subseteq K_t.
\end{equation}
Moreover, we  have the following.

\noindent \textbf{Assertion B.}  \ There exists a subset $M\subseteq [0,t-1]$ with $|M|\geq t-\delta$ such that $K_i\subsetneq K_{i+1}\mbox{ for each } i\in M.$

{\sl Proof of Assertion B.} \ Let $v\in [0,t-1]$ be an arbitrary index with $K_v=K_{v+1}$. We shall apply Lemma \ref{Lemma three conclusions} by taking $a=\sum\limits_{i=1}^{v+1} a_i$ and $b=\sum\limits_{i=1}^{v} a_i$. Since $a<_{\mathcal{H}} b$, it follows from Lemma \ref{Lemma three conclusions} that $\gcd(\kappa_j(\theta_b),n_j)\mid  \gcd(\kappa_j(\theta_a),n_j)$ for each $j\in [1,r]$, and moveover, $${\rm pot}_2(\gcd(\kappa_w(\theta_b), \ n_w))<{\rm pot}_2(\gcd(\kappa_w(\theta_a), \ n_w))={\rm pot}_2(n_w)$$ for some $w\in [1,r]$.
By the arbitrariness of $v$, we have Assertion B proved. \qed

For each $m\in M$, since $\frac{{\rm U}(\mathcal{S}_R)}{K_{m+1}}\cong \frac{{\rm U}(\mathcal{S}_R)/ K_{m}}{K_{m+1}/ K_{m}}$ and ${\rm D}(K_{m+1}/ K_{m})\geq 2$, it follows from Lemma \ref{Lemma recusive Davenport constant} that
\begin{equation}\label{equation D()leq D()-1}
\begin{array}{llll}
{\rm D}({\rm U}(\mathcal{S}_R)/ K_{m+1})&=& {\rm D}(\frac{{\rm U}(\mathcal{S}_R)/ K_{m}}{K_{m+1}/ K_{m}}) \\
&\leq & {\rm D}({\rm U}(\mathcal{S}_R)/ K_{m})-({\rm D}(K_{m+1}/ K_{m})-1)\\
&\leq & {\rm D}({\rm U}(\mathcal{S}_R)/ K_{m})-1.\\
\end{array}
\end{equation}
By \eqref{equation subgroup chains},  \eqref{equation D()leq D()-1} and Assertion B,  we conclude that \begin{align}\label{equation length and D()}
\begin{array}{llll}
1\leq {\rm D}({\rm U}(\mathcal{S}_R)/ K_t)&\leq & {\rm D}({\rm U}(\mathcal{S}_R)/ K_{0})-|M| \\
&\leq & {\rm D}({\rm U}(\mathcal{S}_R))-(t-\delta)\\
&=& (\ell-\delta)-(t-\delta)\\
&=& \ell-t\\
&=& |TV^{[-1]}|.\\
\end{array}
\end{align}
By \eqref{equation sigma(V)Hsigma(T)} and Lemma \ref{Lemma three conclusions} {\rm (ii)}, we have
\begin{equation}\label{equation two common divisors equal}
\gcd(\kappa_i(\theta_{\sigma(V)}),n_i)=\gcd(\kappa_i(\theta_{\sigma(T)}),n_i)
\end{equation} for each $i\in [1,r]$.
Let $$\mathcal{P}_i=\{p: p \mbox{ is prime with } {\rm pot}_p(\gcd(\kappa_i(\theta_{\sigma(V)})), n_i)={\rm pot}_p(n_i)>0\}$$ where $i\in [1,r]$.

Let $a$ be an arbitrary term of $TV^{[-1]}$. By \eqref{equation two common divisors equal}, we have that for each $i\in [1,r]$,
\begin{equation}\label{equation q notmid a}
q\not\mid\kappa_i(\theta_{a})
\end{equation}
for any prime divisor $q$ of $n_i$ with $q\notin\mathcal{P}_i$. By the Chinese Remainder Theorem, we can choose integers $\tilde{a}_i\in [0,n_i-1]$ such that
\begin{equation}\label{equation tilde a 1}
\tilde{a}_i\equiv 1 \pmod {p^{{\rm pot}_p(n_i)}} \ \ \ \mbox{ for any prime } p\in \mathcal{P}_i
\end{equation}
and
\begin{equation}\label{equation tilde a 2}
\tilde{a}_i\equiv \kappa_i(\theta_{a}) \pmod {q^{{\rm pot}_q(n_i)}} \ \ \ \mbox{ for any prime divisor } q \mbox{ of } n_i  \mbox{ with } q\notin \mathcal{P}_i,
\end{equation}
where $i\in [1,r]$.
Let $\tilde{a}$ be the element of $\mathcal{S}_R$ with $\kappa_i(\theta_{\tilde{a}})=\tilde{a}_i \mbox{ for each }i\in [1,r].$
By \eqref{equation q notmid a}, \eqref{equation tilde a 1} and \eqref{equation tilde a 2}, we conclude that $\gcd(\kappa_i(\theta_{\tilde{a}}), n_i)=1$ for each $i\in [1,r]$, i.e.,
\begin{equation}\label{equation tilde a in U}
\tilde{a}\in {\rm U}(\mathcal{S}_R).
\end{equation}
By  \eqref{equation tilde a 1} and \eqref{equation tilde a 2}, we conclude that
$$\tilde{a}_i*\kappa_i(\theta_{\sigma(V)})\equiv 1\cdot \kappa_i(\theta_{\sigma(V)})\equiv 0\equiv \kappa_i(\theta_{a})*0\equiv \kappa_i(\theta_{a})*\kappa_i(\theta_{\sigma(V)})\pmod {p^{{\rm pot}_p(n_i)}}$$ for any prime $p\in \mathcal{P}_i,$
and that
$$\tilde{a}_i*\kappa_i(\theta_{\sigma(V)})\equiv \kappa_i(\theta_{a})* \kappa_i(\theta_{\sigma(V)})\pmod {q^{{\rm pot}_q(n_i)}}$$ for any prime divisor $q$ of $n_i$ with $q\notin \mathcal{P}_i,$ that is,
\begin{equation}\label{equation sigma(V)+tilde a}
\sigma(V)+\tilde{a}=\sigma(V)+a \ \ \ \mbox{for each term}\ \  a \ \mbox{of} \ TV^{[-1]}.
\end{equation}

By \eqref{equation length and D()}, \eqref{equation tilde a in U} and the arbitrariness of the element $a$ above, we see that $\coprod\limits_{a\mid TV^{[-1]}}\tilde{a}$ is a nonempty sequence of terms from ${\rm U}(\mathcal{S}_R)$ of length $|\coprod\limits_{a\mid TV^{[-1]}}\tilde{a}|=|TV^{[-1]}|\geq {\rm D}({\rm U}(\mathcal{S}_R)/ K_t)$. It follows that there exists a {\sl nonempty} subsequence $W\mid TV^{[-1]}$  such that
$\sigma(\coprod\limits_{a\mid W}\tilde{a})\in K_t$ which implies
\begin{equation}\label{equation sigma(tilde a)=0}
\sigma(V)+\sigma(\coprod\limits_{a\mid W}\tilde{a})=\sigma(V).
\end{equation}
By \eqref{equation sigma(V)+tilde a} and \eqref{equation sigma(tilde a)=0}, we conclude that
$$\begin{array}{llll}
\sigma(T)&=& \sigma(TW^{[-1]}V^{[-1]})+(\sigma(V)+\sigma(W)) \\
&=& \sigma(TW^{[-1]}V^{[-1]})+(\sigma(V)+\sigma(\coprod\limits_{a\mid W}\tilde{a}))\\
&=& \sigma(TW^{[-1]}V^{[-1]})+\sigma(V)\\
&=& \sigma(TW^{[-1]}),\\
\end{array}$$
and $T'=TW^{[-1]}$ is the desired proper subsequence of $T$.
This completes the proof of the theorem. \qed

\section{Concluding remarks}

We remark that the upper bound in Theorem \ref{Theorem main} could be reached in some cases. For example, take $$R=\mathbb{Z}_8^{r_1}\oplus\mathbb{Z}_4^{r_2}\oplus\mathbb{Z}_2^{r_3}.$$
Note that ${\rm U}(\mathcal{S}_R)\cong C_2^{2r_1+r_2}$ is the direct sum of $2r_1+r_2$ cyclic groups of order two.  It is well known that ${\rm D}({\rm U}(\mathcal{S}_R))=2r_1+r_2+1.$ For each $i\in [1,r_1+r_2+r_3]$, we take $a_i\in \mathcal{S}_R$ with $\theta_{a_i}=(1,\ldots,1,2,1\ldots,1)$ where the only $2$ occurs at the $i$-th location. It is not hard to check that
the sequence $$T=(\coprod\limits_{i=1}^{r_1}a_i^{[3]})\cdot (\coprod\limits_{j=r_1+1}^{r_1+r_2}a_j^{[2]})\cdot (\coprod\limits_{k=r_1+r_2+1}^{r_1+r_2+r_3}a_k)$$ is an irreducible sequence in the semigroup $\mathcal{S}_R$. It follows that ${\rm D}(\mathcal{S}_R)\geq |T|+1=3r_1+2r_2+r_3+1={\rm D}({\rm U}(\mathcal{S}_R))+r_1+r_2+r_3= {\rm D}({\rm U}(\mathcal{S}_R))+\delta,$ and thus, $${\rm D}(\mathcal{S}_R)={\rm D}({\rm U}(\mathcal{S}_R))+\delta.$$

However, it would be attractive to determine the precise value of
${\rm D}(\mathcal{S}_R)-{\rm D}({\rm U}(\mathcal{S}_R))$ when $R=\mathbb{Z}_{n_1}\oplus\cdots\oplus \mathbb{Z}_{n_r}$
for any integers $r\geq 1$, $n_1,\ldots,n_r>1$.
Even, the following conjecture seems to be interesting.

\begin{conj} \  For integers $r\geq 1$, $n_1,\ldots,n_r>1$, let $R=\mathbb{Z}_{n_1}\oplus\cdots\oplus \mathbb{Z}_{n_r}$. Then $${\rm D}(\mathcal{S}_R)-{\rm D}({\rm U}(\mathcal{S}_R))\leq \sharp \{i\in [1,r]: {\rm pot}_2(n_i)\in [1,3]\}.$$
\end{conj}

\noindent {\bf Acknowledgements}

\noindent
The authors thank Calvin Deng in 2015 for writing us to tell us the ring $\mathbb{Z}_4$ is a counterexample for our Theorem $D$. Before this, in 2010 in  the process of dealing with the generalization of Erd\H{o}s-Ginzburg-Ziv theorem for finite commutative semigroups with Prof. Adhikari,  the two authors have noticed such mistake. As first, we wish to give a precise result of the Davenport constant, however, we failed. Recently, more additive research on semigroups have been done, the authors don't hesitate to stop this error and to propose the remaining problem.

This work is supported by NSFC (11301381, 11271207), Science and Technology Development Fund of Tianjin Higher
Institutions (20121003).


\begin{thebibliography}{99}


\bibitem{AdhikariGaoWang14} S.D. Adhikari, W.D. Gao and G.Q. Wang, \emph{Erd\H{o}s-Ginzburg-Ziv theorem for finite commutative semigroups,} Semigroup Forum, \textbf{88} (2014)  555--568.



\bibitem{GH} A. Geroldinger and F. Halter-Koch, \emph{Non-Unique
Factorizations. Algebraic, Combinatorial and Analytic Theory,}
Pure and Applied Mathematics, vol. 278, Chapman $\&$ Hall/CRC,
2006.

\bibitem{rog1}
K. Rogers, \emph{A Combinatorial problem in Abelian groups,} Proc.
Cambridge Phil. Soc., \textbf{59} (1963) 559--562.

\bibitem{wangDavenportII}  G.Q. Wang, \emph{Davenport constant for semigroups II,}  J. Number Theory, \textbf{153} (2015) 124--134.

\bibitem{wangErdosBugess}  G.Q. Wang, \emph{Structure of the largest idempotent-free sequences in finite semigroups,}  arXiv:1405.6278.

\bibitem{wangAddtiveirreducible}  G.Q. Wang, \emph{Additively irreducible sequences in commutative semigroups,}  arXiv:1504.06818.

\bibitem{wanggao} G.Q. Wang and W.D. Gao,
\emph{Davenport constant for semigroups,} Semigroup Forum,
\textbf{76} (2008) 234--238.


\bibitem{wang-zhang-wang-qu}  H.L. Wang, L.Z. Zhang, Q.H. Wang and Y.K. Qu,  \emph{Davenport constant of the multiplicative semigroup of the quotient ring $\frac{\F_p[x]}{\langle f(x)\rangle}$,} International Journal of Number Theory, in press,
 DOI: 10.1142/S1793042116500433.


 \bibitem{wang-zhang-qu} L.Z. Zhang, H.L. Wang and Y.K. Qu,  \emph{A problem of Wang on Davenport constant for the multiplicative semigroup of the quotient ring of $\F_2[x]$},  arXiv:1507.03182.

\end{thebibliography}
\end{document}